\documentclass[reqno,12pt]{amsart}
\usepackage{amsmath,amssymb,amscd,amsxtra,dsfont}
\usepackage{graphicx,pst-grad,pst-plot,pst-coil}
\usepackage{pstricks}
\usepackage{eucal}
\usepackage{mathtools}
\usepackage[margin=1.3in]{geometry}

\newtheorem{theorem}{Theorem}
\newtheorem{proposition}[theorem]{Proposition}
\newtheorem{cor}[theorem]{Corollary}
\newtheorem{lemma}[theorem]{Lemma}

\theoremstyle{remark}

\usepackage{mathptmx}

\usepackage[all]{xy}
\def\bfig{\vcenter\bgroup}
\def\efig{\egroup}

\def \defi{\begin{definition}}
\def \edefi{\end{definition}}

\def \prop{\begin{proposition}}
\def \eprop{\end{proposition}}

\def \thm{\begin{theorem}}
\def \ethm{\end{theorem}}

\def \lem{\begin{lemma}}
\def \elem{\end{lemma}}

\def \cor{\begin{corollary}}
\def \ecor{\end{corollary}}

\newenvironment{definition}{\refstepcounter{theorem}
{\medskip\par\noindent\bf Definition \arabic{section}.\arabic{theorem}.
}}{\vskip 2ex\par}
\newenvironment{example}{\refstepcounter{theorem}
{\medskip\par\noindent\it Example \arabic{section}.\arabic{theorem}.}}{\vskip 2ex\par}

\newenvironment{example*}{
{\medskip\noindent\it Example.}}{\vskip 2ex\par}

\newcommand{\ex}{\begin{example}}
\newcommand{\eex}{\end{example}}

\newenvironment{exer*}
  {\small\begin{exercise}}
  {\end{exercise}}

\def \ex*{\begin{example*}}
\def \eex*{\end{example*}}

\newenvironment{remark*}{
{\medskip\noindent\it Remark.}}{\vskip 2ex\par}
\def \rem*{\begin{remark*}}
\def \erem*{\end{remark*}}

\newenvironment{claim*}{
{\medskip\noindent\it Claim.}}{\vskip 2ex\par}

\def \pf{\begin{proof}}
\def \epf{\end{proof}}

\def \enum{\begin{enumerate}}
\def \eenum{\end{enumerate}}

\numberwithin{equation}{section}
\numberwithin{figure}{section}

\newcommand{\C}{\mathbb{C}}

\newcommand{\g}{\mathfrak{g}}

\newcommand{\Z}{\mathbb{Z}}
\newcommand{\ga}{\alpha}

\newcommand{\gc}{\gamma}

\newcommand{\R}{\mathbb{R}}

\newcommand{\ddtz}{\left. \frac{d{ }}{dt}\right|_{t=0} }

\def\dual{^{\scriptstyle\vee}}

\newcommand{\Ad}{\operatorname{Ad}}

\newcommand{\cinf}{{C^{\infty}}}

\newcommand{\Fl}{\operatorname{F\ell}}

\newcommand{\inv}{^{-1}}

\newcommand{\pt}{\operatorname{pt}}

\newcommand{\term}[1]{\textbf{\textit{#1}}}

\newcommand{\Unit}{\text{U}}

\begin{document}
\title
{Computing Topological Invariants Using Fixed Points
}
\author{Loring W. Tu}
\address{Department of Mathematics\\
         Tufts University\\
         Medford, MA 02155}
\email{loring.tu@tufts.edu}
\date{April 17, 2016; version 10.  This article is based on a talk given at the Sixth
International Congress of Chinese Mathematicians, Taipei, Taiwan, in 2013}
\begin{abstract}
When a torus acts on a compact oriented manifold with isolated fixed points,
the equivariant localization formula of Atiyah--Bott and Berline--Vergne converts
the integral of an equivariantly closed form into a finite sum over the fixed points
of the action, thus providing a powerful tool for computing integrals on a manifold.
An integral can also be viewed as a pushforward map from a manifold to a point,
and in this guise it is intimately related to the Gysin homomorphism.
This article highlights two applications of the equivariant localization formula.
We show how to use it to compute characteristic numbers of a homogeneous space
and to derive a formula for the Gysin map of a fiber bundle.

\end{abstract}
\maketitle

Many invariants in geometry and topology can be represented as integrals.
For example, according to the Gauss--Bonnet theorem, the Euler characteristic
of a compact oriented surface in $\R^3$ is $1/2\pi$ times the integral of its
Gaussian curvature:
\[
\chi (M) = \frac{1}{2\pi} \int_M K {\rm vol}.
\]

The Euler characteristic can be generalized to other characteristic numbers.
For example, if $E$ is a complex vector bundle of rank $r$ over a
complex manifold $M$ of complex dimension $n$, and $c_1, \ldots, c_r$ are the Chern classes of $E$, then
the integrals
\[
\int_M  c_1^{i_1} \cdots c_r^{i_r}, \text{ where } \sum_{k=1}^r  k \cdot i_k = n,
\]
are the \term{Chern numbers} of $E$.
Taking $E$ to be the holomorphic tangent bundle $TM$ of $M$,
the Chern numbers of $TM$ are called the \term{Chern numbers} of the complex manifold $M$.
They are smooth invariants of $M$.
The top Chern number $\int_M c_n(TM)$ is the Euler characteristic $\chi(M)$.

In general, integrals on a manifold are notoriously difficult to compute, but if there is a 
torus action on the manifold with isolated fixed points, then the equivariant localization formula of Atiyah--Bott and Berline--Vergne converts
certain integrals into finite sums over the fixed points.

A \term{homogeneous space} is a space of the form $G/H$, where
$G$ is a Lie group and $H$ is a closed subgroup.
We will first consider the following problem:

\bigskip

\noindent
\textbf{Problem 1}.  How does one compute an integral on a homogeneous
space $G/H$?

\bigskip

Every Lie group has a maximal torus. 
Since the maximal tori in a
Lie group are all conjugate to one another,
they all have the same dimension.
The dimension of a maximal torus in a Lie group $G$
is called the \term{rank} of $G$.
The method outlined in this article applies when $G$
is a compact connected Lie group and $H$ is a closed
subgroup of maximal rank.
In this case, a maximal torus of $H$ is also a maximal
torus of $G$.
The complex projective spaces $\C P^n$, 
the complex Grassmannians  $G(k, \C^n)$,
and complex flag manifolds are all examples of such homogeneous spaces.

To simplify the exposition, we will assume throughout
that homology and cohomology are taken with real
coefficients.
Cohomology with real coefficients is to be interpreted as
singular cohomology or de Rham cohomology, as the case may be.
Suppose $f\colon E \to M$ is a continuous map between compact
oriented manifolds of dimensions $e$ and $m$ respectively.
Then there is an induced map $f_*\colon H_*(E) \to
H_*(M)$ in homology, and by Poincar\'e duality, an induced map
$H^{e-*}(E) \to H^{m-*}(M)$ in cohomology.
This map in cohomology, also denoted by $f_*$, is called the
\term{Gysin map}.
It is defined by the commutative diagram
\[
\xymatrix{
H^k(E) \ar[r]^-{f_*} \ar[d]_-{{\rm P.D.}}^-{\simeq} & 
H^{k-(e-m)}(M) \ar[d]^-{{\rm P.D.}}_-{\simeq} \\
H_{e-k}(E) \ar[r]_-{f_*} & H_{e-k}(M),
}
\]
When $M$ is a point, the Gysin map $f_*\colon H^e(E) \to H^0(\text{pt}) = \R$
is simply integration.  Thus, the Gysin map generalizes integration.
For a fiber bundle $f\colon E \to M$, the Gysin map $f_*\colon H^k(E)
\to H^{k-(e-m)}(M)$ is integration along the fiber; it lowers the
degree
by the fiber dimension $e-m$.

\bigskip

\noindent
\textbf{Problem 2}.  Derive a formula for the Gysin map of a fiber bundle
$f\colon E \to M$ whose fibers are homogeneous spaces.

\bigskip

In fact, our method solves Problem 2 for other fiber bundles as well;
it suffices that the fibers be \term{equivariantly formal},
as defined below.

Because this article is meant to be expository, we do not give
complete proofs of most of the results cited.
For further details, proofs, and references, consult \cite{tu10}
and \cite{tu13}.

\bigskip

\noindent
{\sc Acknowledgments}.  The author is grateful to Jeffrey D.\ Carlson
for his careful reading of several drafts of the manuscript and for his
helpful comments, to the Tufts Faculty Research Award Committee
for its financial support, and to
the National Center for Theoretical Sciences Mathematics Division (Taipei Office) in Taiwan for hosting him
during part of the preparation of this manuscript.

\section{The Fixed Points of a Torus Action on $G/T$}

To apply the equivariant localization formula, we need a torus action
on a manifold.  In fact, any action by a compact Lie group will do,
because every compact Lie group contains a maximal torus,
and we can simply restrict the given action to that of a 
maximal torus.

For simplicity, I will assume that in the homogeneous space $G/H$,
the closed subgroup $H$ is a maximal torus $T$ in the Lie group $G$.
The method of treating the general case is similar, but the
formulas are a bit more complicated.

The torus $T$ acts on $G/T$ by left multiplication:  
\begin{align*}
T \times G/T &\to G/T,\\
t \cdot xT &= txT.
\end{align*}
Let us compute the fixed point set of this action:
\begin{align*}
\text{$xT$ is a fixed point} \quad 
\Longleftrightarrow &\quad t\cdot xT = xT \ \text{ for all } t \in T\\
\Longleftrightarrow &\quad x\inv t x T = T \ \text{ for all } t \in T\\
\Longleftrightarrow &\quad x\inv T x \subset T \\
\Longleftrightarrow &\quad x \in N_G(T) = \text{normalizer of $T$ in $G$}\\
\Longleftrightarrow &\quad xT \in N_G(T)/T.
\end{align*}
The group $W\coloneqq N_G(T)/T$ is called the \term{Weyl group} of
$T$ in $G$ and is well known from the theory of Lie groups to be a finite
reflection group.
Thus, the action of $T$ on $G/T$ by left multiplication has finitely many
fixed points.

\section{Equivariant Cohomology}

To study the algebraic topology of spaces with group actions, 
one looks for a functor that incorporates in it both the topology
of the space and the action of the group.
Let $M$ be a topological space on which a topological group $G$ acts.
Such a space is called a \term{$G$-space}.
Equivariant cohomology is a functor from the category of
$G$-spaces to the category of commutative rings.

As a first candidate, one might try the singular cohomology $H^*(M/G)$
of the orbit space $M/G$.
Consider the following two examples.

\begin{example}
The circle $G=S^1$ acts on the $2$-sphere $M = S^2$ in $\R^3$
by rotation about the $z$-axis.
Each orbit is a horizontal circle and the orbit space $M/G$
is homeomorphic to the closed interval $[-1,1]$ on the $z$-axis.
The cohomology $H^*(M/G)$ is trivial.
\end{example}

\begin{example}
The group $G= \Z$ of integers acts on $M= \R$ by
translations  $n \cdot x = x+n$.
The orbit space $M/G$ is the circle $\R/\Z$.
The cohomology $H^*(M/G)$ of the orbit space is
nontrivial.
\end{example}

In the first example, the cohomology of the orbit space $M/G$
yields little information about the action.
In the second example, $H^*(M/G)$ provides an interesting
invariant of the action.
An action of $G$ on a space $M$ is said to be
\term{free} if the stabilizer of every point is the
identity.  One difference between the two examples is
that the action of $S^1$ on $S^2$ by rotation is not
free---the stabilizers at the north and south poles are the 
group $S^1$ itself---while the action of $\Z$ on $\R$ by
translation is free.
In general, when the action of $G$ on $M$ is not free,
the quotient space $M/G$ may be problematical.

In homotopy theory there is a standard procedure for 
converting a nonfree action to a free action.
Suppose $EG$ is a contractible space on which $G$
acts freely.
Then no matter how $G$ acts on the space $M$,
the diagonal action of $G$ on $EG \times M$ will
be free:
\begin{align*}
(e,x) = g \cdot (e,x) = (g\cdot e, g\cdot x) &\quad\Longleftrightarrow\quad  g\cdot e = e \text{ and } g\cdot x = x\\
&\quad\Longrightarrow\quad g\cdot e = e \\
&\quad\Longrightarrow\quad g= 1.
\end{align*}
Since $EG$ is contractible, $EG \times M$ has the
same homotopy type as $M$.
Homotopy theorists call the orbit space of $EG \times M$
under the free diagonal action of $G$ the
\term{homotopy quotient} $M_G$ of $M$ by $G$.
The \term{equivariant cohomology} $H_G^*(M)$ is
defined to be the singular cohomology $H^*(M_G)$
of the homotopy quotient, with whatever coefficient ring
is desired.
It can be shown that the homotopy type of $M_G$ is
independent of the choice of $EG$.
Thus, the equivariant cohomology $H_G^*(M)$ is well defined.

From the theory of principal bundles, we know that a
weakly contractible space $EG$ on which a topological group $G$
acts freely is the total
space of a universal principal $G$-bundle
$EG \to BG$, a principal $G$-bundle from which
any principal $G$-bundle can be pulled back.
That is, given any principal $G$-bundle $P \to M$,
there is a map $f\colon M \to BG$ such that $P$ is
isomorphic to the pullback bundle $f^*EG$.

The process of constructing from a universal 
principal $G$-bundle $\alpha\colon EG \to BG$
and a left $G$-space $M$ the 
homotopy quotient $M_G = (EG \times M)/G$ is called
the \term{Borel construction}.
We denote by $[e,x]$ the equivalence class in $M_G$ of
$(e,x) \in EG \times M$.
It is easy to check that the natural map $M_G \to BG$, $[e,x] \mapsto \alpha(e)$,
is a fiber bundle with fiber $M$ and structure group $G$.
The inclusion of the fiber $M$ into $M_G$ induces a homomorphism
$H^*(M_G) \to H^*(M)$ in cohomology.
Hence, there is a canonical map $H_G^*(M) \to H^*(M)$
from equivariant cohomology to ordinary cohomology.
A cohomology class in $H^*(M)$ in the image of this map is
said to have an \term{equivariant extension}.

A vector bundle $\pi\colon V \to M$ is said to be
\term{$G$-equivariant} if $V$ and $M$ are left $G$-spaces
and $\pi\colon V \to M$ is a $G$-map such that for
every $g \in G$ and every fiber $V_x$, the map
$\ell_g\colon V_x \to V_{gx}$ is a linear map.

There is one situation in which a cohomology class on $M$
automatically has an equivariant extension, namely
when it is a characteristic class $c(V)$ of a $G$-equivariant
vector bundle $V \to M$.
In this case, the induced map $V_G \to M_G$ is a vector
bundle with the same fiber as $V \to M$ and the equivariant
characteristic class $c^G(V)$ is defined to be
\[
c^G(V) \coloneqq c(V_G) \in H_G^*(M).
\]
The commutative diagram
\[
\xymatrix{
V \ar@{^{(}->}[r] \ar[d] & V_G \ar[d]\\
M \ar@{^{(}->}[r]_j & M_G}
\]
shows that the bundle $V \to M$ is the restriction of
$V_G \to M_G$ to $M$; i.e.,
if $j\colon M \to M_G$ is the inclusion,
then $V = j^*(V_G)$.
By the naturality of characteristic classes,
\[
c(V) = c(j^* V_G)= j^* c(V_G) = j^* c^G(V).
\]
Thus, the cohomology class $c(V)$ has equivariant 
extension $c^G(V)$.

In general, the canonical map $H_G^*(M) \to H^*(M)$ is neither
surjective nor injective.
If it is surjective, then $M$ is said to be \term{$G$-equivariantly
formal}.
A $G$-equivariantly formal space $M$ is then one in which every
cohomology class in $H^*(M)$ has an equivariant extension.
It turns out that any
homogeneous space $G/H$, where $H$ contains a maximal torus of $G$,
is equivariantly formal.

\begin{example}
Any group $G$ acts trivially on a point.
The homotopy quotient of a point by $G$ is
\[
\pt_G = EG \times_G \pt = EG/G = BG,
\]
so the equivariant cohomology of a point is $H^*(BG)$.
\end{example}

For any $G$-space $M$, the constant map $M \to \pt$
induces a ring homomorphism $H_G^*(\pt) \to H_G^*(M)$,
which shows that the equivariant cohomology ring $H_G^*(M)$ has
the structure of an algebra over the ring $H^*(BG)$.

\begin{example}
The circle $S^1$ acts freely on the sphere $S^{2n+1}$ with
quotient $\C P^n$. 
Therefore, it acts
freely on the union $S^{\infty} = \bigcup_{n=1}^{\infty} S^{2n+1}$ 
with quotient $\C P^{\infty} = \bigcup_{n=1}^{\infty} \C P^n$.
Since the homotopy groups $\pi_k(S^{\infty})$ vanish for all $k$,
by Whitehead's theorem, $S^{\infty}$ is contractible.
Thus, $ES^1 = S^{\infty}$ and $BS^1 = \C P^{\infty}$, and
\[
H_{S^1}^*(\pt) = H^*(BS^1) = H^*(\C P^{\infty}) = \R [u].
\]
\end{example}

\section{Equivariant Forms}

Just as the singular cohomology with real coefficients
of a manifold can be computed using differential forms
(\term{de Rham's theorem}),
the equivariant cohomology of a manifold with a 
group action can be computed using equivariant differential forms
(\term{Cartan's theorem}).

Suppose a Lie group $G$ acts on a manifold $M$.
Fix a basis $B_1, \ldots, B_m$ for the Lie algebra $\g$,
and let $v_1, \ldots, v_m$ be the dual basis for $\g\dual{}$.
A function $\ga\colon \g \to \Omega(M)$ is said to be \term{polynomial}
if it can be written as a polynomial in $v_1, \ldots, v_m$ with
coefficients that are $\cinf$ forms on $M$:
\[
\ga = \sum \ga_I v_1^{i_1} \cdots v_m^{i_m}, \quad \ga_I \in \Omega(M).
\]
For $A \in \g$, the value of $\ga$ at $A$ is given by
\[
\ga(A) = \sum \ga_I v_1(A)^{i_1} \cdots v_m(A)^{i_m}, \quad \ga_I \in \Omega(M).
\]

The group $G$ acts on $\g$ by the adjoint representation
and on $\Omega(M)$ by the pullback: $g\cdot \omega = \ell_{g^{-1}}^*\omega$.
A \term{$G$-equivariant form} on $M$ is a polynomial
map $\ga\colon \g \to \Omega(M)$ $G$-equivariant with respect to 
these actions:
\[
\ga \big( (\Ad g) A\big) = \ell_{g\inv}^* \ga(A)
\]
for all $g\in G$ and $A \in \g$.
If $G$ is a torus $T$, then the adjoint action is trivial and
a $T$-equivariant form is a polynomial in $v_1, \ldots, v_m$
with $T$-invariant forms on $M$ as coefficients:
\[
\ga = \sum \ga_I v_1^{i_1} \cdots v_m^{i_m},  \quad \ga_I \in \Omega(M)^T.
\]

Each element $A \in \g$ gives rise to a vector field $\underline{A}$
on
the $G$-manifold $M$, called a \term{fundamental vector field},
by
\[
\underline{A}_p = \ddtz e^{-tA} \cdot p \quad \text{for $p \in M$}.
\]
Then $\g$ acts on $\Omega(M)$ by $\iota_A \omega \coloneqq \iota_{\underline{A}}\omega$
for $A \in \g$ and $\omega \in \Omega(M)$.
Let $\Omega_G(M)$ be the set of $G$-equivariant forms on $M$.
It is an algebra over $\R$ equipped with an operator $d_G$ whose square
is zero, called the
\term{Cartan differential}, given by
\[
(d_G \ga)(A) = d\big(\ga(A)\big) - \iota_{A} \big(
\ga(A)\big)\quad \text{for all } A \in \g.
\]
In terms of the basis $B_1, \ldots, B_m$ for $\g$
and the dual basis $v_1, \ldots, v_m$ for $\g\dual$,
we may write $A = \sum v_i(A) B_i$ and
\[
(d_G\ga)(A) = \sum v_1(A)^{i_1} \cdots v_m(A)^{i_m} d\ga_I - 
\sum v_i(A)\iota_{B_i} v_1(A)^{i_1} \cdots v_m(A)^{i_m}  \ga_I.
\]
Hence,
\begin{align*}
d_G\ga &= \sum v_1^{i_1} \cdots v_m^{i_m} d\ga_I - \sum v_i\iota_{B_i} v_1^{i_1} \cdots v_m^{i_m} \ga_I\\
&= d\big( \sum v^I \ga_I\big) - \sum  v_i\iota_{B_i}\big(\sum v^I {\ga_I}\big)\\
&= d \ga - \sum v_i \iota_{B_i} \ga.
\end{align*}

To integrate an equivariant form, one simply integrates its coefficients:
\[
\int_M \sum \ga_I v^I = \sum \Big(\int_M \ga_I\Big) v^I.
\]
The integrals $\int_M \omega$ on a manifold $M$ of dimension $n$ that are amenable to
computation using the equivariant localization formula are integrals of differential
$n$-forms that have equivariantly closed extensions $\tilde{\omega}$:
\[
\tilde{\omega} = \omega + \sum \ga_I v^I,\quad d_G\tilde{\omega} = 0.
\]
For dimension reasons,
\[
\deg \ga_I = n- 2 \sum |I| < n,
\]
so that 
\[
\int_M \tilde{\omega} 
= \int_M \omega + \sum \Big( \int \ga_I\Big) v^I = \int_M \omega.
\]
The equivariant localization formula then gives the integral $\int_M \tilde{\omega}$ as a finite sum.

\section{Line Bundles on $G/T$ and $BT$}
The quotient map $G \to G/T$ is a principal $T$-bundle.
Recall that every principal $T$-bundle is a pullback from the universal $T$-bundle
$ET \to BT$.
A \term{character} of a torus $T$ is a homomorphism $\gamma\colon T \to \C^{\times}$,
which can be viewed as an action of $T$ on $\C$.
By the mixing construction, we can associate to a character $\gamma$
a complex line bundle $L_{\gamma} \coloneqq  G \times_{\gamma} \C$ on $G/T$
and a complex line bundle $S_{\gamma} \coloneqq  ET \times_{\gamma} \C$ on
$BT$.
Since $G \to G/T$ is a pullback of $ET \to BT$,
the line bundle $L_{\gc}$ is the pullback of $S_{\gc}$ by
the classifying map $G/T \to BT$.
Left multiplication by elements of $T$ makes $L_{\gc}$ into a $T$-equivariant
complex line bundle over $G/T$.

\section{Cohomology Classes on $G/T$ and $BT$}
The first Chern class $c_1(L_{\gamma})$ is a cohomology class of degree 2
on $G/T$.  
Choose a basis $\chi_1, \ldots, \chi_{\ell}$ for the characters of $T$ and let
\[
y_i = c_1 (L_{{\chi}_i})\in H^2(G/T), \qquad u_i  = c_1 (S_{\chi_i}) \in H^2(BT).
\]

The computation of an integral over $G/T$ using the equivariant localization
formula is made possible by the happy fact that the cohomology ring of $G/T$
is generated by the $y_i$.
As characteristic classes of $T$-equivariant bundles, these Chern classes automatically
have equivariantly closed extensions, namely, the equivariant Chern classes, 
so the integrals of monomials in the Chern classes
 can be calculated using the equivariant localization formula.

If $T = \overbrace{S^1 \times \cdots \times S^1}^{\ell \text{
    times}}$, then its classifying space is
\begin{align*}
BT &=  BS^1 \times \cdots \times BS^1 \\
&= \C P^{\infty} \times \cdots \times \C P^{\infty}.
\end{align*}
By the K\"unneth formula, the cohomology ring of $BT$ is
\[
H^*(BT) = \R[u_1, \ldots, u_{\ell}].
\]

\section{The Action of the Weyl Group on the Polynomial Ring $H^*(BT)$} \label{6s:action}
Let $N_G(T)$ be the normalizer of $T$ in $G$.
Recall that the Weyl group of $T$ in $G$ is
\[
W : = W_G(T) = N_G(T) /T.
\]
Given a principal $T$-bundle $X \to X/T$, there is always an
action of the Weyl group $W$ on the base $X/T$ by
\[
(xT) w = xwT.
\]
Hence, $W$ acts on the base space $BT = ET/T$ of the universal bundle $ET$
and there is an induced action on the polynomial ring
$H^*(BT)=\R[u_1, \ldots, u_{\ell}]$.

\section{The Cohomology Ring of $G/T$}
Let $R$ be the polynomial ring 
\[
R = \R[y_1, \ldots, y_{\ell}] \simeq \R[u_1, \ldots, u_{\ell}].
\]
The Weyl group $W$ acts on $R$, as explained in Section~\ref{6s:action}.
Let $(R_+^W)$ be the ideal generated by the invariant
polynomials of positive degree.

\thm[{\cite[Prop.\ 26.1, p.\ 190]{borel53}, \cite[Th.~5, p.~190]{tu10}}]
 The cohomology ring of $G/T$ is
\[
H^*(G/T) = \frac{R}{(R_+^W)} = \frac{\R[y_1, \ldots, y_{\ell}]}
{(\R[y_1, \ldots, y_{\ell}]_+^W)}.
\]
\ethm

\section{The Equivariant Localization Formula}
Suppose a torus acts smoothly on a manifold $M$ with
isolated fixed points and $\tilde{\omega}$ is an
equivariantly closed form.  Let $i_p\colon \{p\} \to M$
be the inclusion of a point. 
It induces a map $(i_p)_T\colon \{p\}_T \to M_T$ of homotopy quotients and
correspondingly a restriction map $(i_p)_T^* \colon H_T^*(M) \to
H_T^*(p)$. 
To simplify the notation, we will write $(i_p)_T^*$ as $i_p^*$.
Then the equivariant localization
formula of Atiyah--Bott \cite{atiyah--bott} and Berline--Vergne \cite{berline--vergne} is the following.

\thm[Equivariant localization formula]
\[
\int_M \tilde{\omega} = \sum_{p \in M^T} 
\frac{i_p^* \tilde{\omega}}{e^T(\nu_p)},
\]
where $\nu_p$ is the normal bundle of $p$ in $M$
and $e^T(\nu_p)$ is the equivariant Euler class of $\nu_p$.
\ethm

\section{Computing Integrals on $G/T$}
Since $H^*(G/T)$ is generated by $y_1, \ldots, y_{\ell}$,
a cohomology class  of degree $n\coloneqq  \dim G/T$
is a homogeneous polynomial 
$f(y_1, \ldots, y_{\ell})$ of degree $n/2$ in the
variables $y_1, \ldots, y_{\ell}$.
Because the $y_i$ are Chern classes of equivariant $T$-bundles,
they all have equivariantly closed extensions $\tilde{y}_i$.
Then
\begin{equation} \label{e1}
\int_{G/T} f(y_1, \ldots, y_{\ell}) = \int_{G/T} f(\tilde{y}_1, \ldots, \tilde{y}_{\ell}),
\end{equation}
which can be computed as a finite sum by the equivariant localization
formula:
\begin{equation} \label{e2}
\int_{G/T} f(\tilde{y}_1, \ldots, \tilde{y}_{\ell})
=\sum_{w \in W} \frac{i_w^* f(\tilde{y})}{e^T(\nu_w)},
\end{equation}
where the sum is taken over the Weyl group $W$ of $T$ in $G$.
To evaluate this sum, it suffices to know the restriction of $f(\tilde{y})$
to a fixed point $w\in W$ and the equivariant Euler class of the normal bundle 
of $w$ in $G/T$.
This is worked out in \cite{tu10}:

\medskip
\noindent
\textbf{Restriction formula} \cite[Prop.~10]{tu10}. 
Let $\chi_1, \ldots, \chi_{\ell}$ be a basis
of the characters of $T$, $L_{\chi_i} = G \times_{\chi_i}\C$ the associated complex line
bundles over $G/T$, and
$S_{\chi_i} = ET \times_{\chi_i} \C$ the associated complex line bundles over the classifying
space $BT$.
If $\tilde{y}_i = c_1^T(L_{\chi_i})$ and $u_i = c_1(S_{\chi_i})$,
then 
\begin{equation} \label{e3}
i_w^* \tilde{y_i} = w \cdot u_i.
\end{equation}

\noindent
\textbf{Euler class formula} \cite[Prop.~13]{tu10}. The equivariant Euler class of the normal
bundle $\nu_w$ at a fixed point $w\in W$ of the left action of $T$ on $G/T$ is
\begin{equation} \label{e4}
e^T(\nu_w)  = w\cdot \bigg(\prod_{\alpha \in \Delta^+} c_1(S_{\alpha})\bigg),
\end{equation}
where $\Delta^+$ is a choice of positive roots for $T$ in $G$.

Putting together \eqref{e1}, \eqref{e2}, \eqref{e3}, \eqref{e4}, we obtain a formula
for an integral over $G/T$ as a finite sum:
\[
\int_{G/T} f(y_1, \ldots, y_{\ell})
= \sum_{w\in W} \frac{w\cdot f(u)}{w\cdot \big(\prod_{\alpha \in \Delta^+} c_1(S_{\alpha})\big)}.
\]

\section{Chern Numbers on a Grassmannian}

Using the same method, we can calculate integrals over $G/H$,
where $G$ is a compact Lie group and $H$ is a closed subgroup
containing a maximal torus of $G$. 
For the complex Grassmannian $G(k, \C^n)$, we find
the following Chern number formula:

\thm
If $S$ is the tautological subbundle over $G(k,\C^n)$, then
\begin{equation} \label{e:charnogkn}
\int_{G(k,\C^n)} c_1(S)^{m_1} \cdots c_k(S)^{m_k} =
\sum_I \dfrac{\prod_{r=1}^k \sigma_r (u_{i_1}, \dots, 
u_{i_k})^{m_r}}
{\prod_{i\in I} \prod_{j\in J} (u_i - u_j)},
\end{equation}
where $\sum m_r = k(n-k)$, $I$ runs over all multi-indices
$1 \le i_1 < \dots < i_k \le n$, $J$ is its complementary multi-index,
and $\sigma_r$ is the $r$th elementary symmetric polynomial.
\ethm

\section{A Chern Number on $\C P^2$}
One of the surprising features of the localization formula 
is that although the right-hand side of \eqref{e:charnogkn} 
is apparently a sum of rational functions of $u_1, \dots, 
u_n$, the sum is in fact an integer.

\textbf{Example}.
As an example, we compute a Chern number on $\C
P^2=G(1,\C^3)$.
The real cohomology of $\C P^2$ is $H^*(\C P^2)= \R[x]/(x^3)$,
generated by $x=c_1(S^{\spcheck}) = -c_1(S)$, where $S$ is the
tautological subbundle on $\C P^2$.
By \eqref{e:charnogkn},
\begin{align*}
\int_{\C P^2} &x^2 = \int_{G(1,\C^3)} c_1(S)^2 = \sum_{i=1}^3
\frac{u_i^2}{\prod_{j\ne i}(u_i - u_j)}\\
=&\frac{u_1^2}{(u_1-u_2)(u_1-u_3)} +\frac{u_2^2}{(u_2-u_1)(u_2-u_3)}
+\frac{u_3^2}{(u_3-u_1)(u_3-u_2)},
\end{align*}
which simplifies to 1, as expected.

\section{Motivation for Studying the Gysin Map}
In enumerative geometry, to count the number of objects satisfying a
set of conditions, one method is to represent the objects satisfying
each condition by cycles in a parameter space $M$ and then to compute
the intersection of these cycles in $M$.
When the parameter space $M$ is a compact oriented manifold,
by Poincar\'e duality, the intersection of cycle classes in homology corresponds
to the cup product in cohomology.
Sometimes, a cycle $B$ in $M$ is the image $f(A)$ of a cycle $A$
in another compact oriented manifold $E$ 
under a map $f\colon E \to M$.
In this case the homology class $[B]$ of $B$ is the image
$f_*[A]$ of the homology class of $A$ under the
induced map $f_*\colon H_*(E) \to H_*(M)$ in homology,
and the Poincar\'e dual $\eta_B$ of $B$ is the image of the Poincar\'e
dual $\eta_A$ of $A$ under the Gysin map.

There are some classical formulas for the Gysin map
of fiber bundles, obtained using various methods depending on 
what the fiber is. 
For example, for a projective bundle $f\colon P(E) \to M$,
where $E \to M$ is a vector bundle, if $\mathcal{O}_{P(E)}(1)$
denotes the dual of the tautological subbundle over $P(E)$,
then
\[
H^*\big(P(E)\big) = H^*(M)[x]/I, 
\]
where $x = c_1 \big( \mathcal{O}_{P(E)}(1)$\big) and $I$ is an ideal in $H^*(M)[x]$.
The formula for the Gysin map 
\cite[Eq.~4.3, p.~318]{arbarello--cornalba--griffiths--harris} is
\[
f_*\bigg(\frac{1}{1-x}\bigg) = \frac{1}{c(E)}.
\]

\section{Homotopy Quotients as Universal Fiber Bundles}
Surprisingly, the equivariant localization formula provides 
a systematic method for computing the Gysin map of a fiber
bundle.
It is based on the fact that the homotopy quotient $F_G$
is the total space of a universal fiber bundle with fiber $F$
and structure group $G$.
Let $E \to M$ be a fiber bundle with fiber $F$ and structure
group $G$.
Associated to $E$ is a principal $G$-bundle $P \to M$ such
that $E$ is the associated bundle $E = P \times_G F$.

The classifying map $\underline{h}$ of the principal bundle
$P\rightarrow M$ in the diagram
\[
\xymatrix{
P \ar[d] \ar[r] & EG \ar[d] \\
M \ar[r]_-{\underline{h}} & BG
}
\]
induces a map of fiber bundles

\begin{equation} \label{e:universal}
\begin{xy}
(2,10.5)*{E =};
(6,5)*{\xymatrix{
P\times_G F \ar[d]_-f \ar[r]^-h & 
EG \times_G F \ar[d]^-{\pi_G} \\
M \ar[r]_-{\underline{h}} & BG.}
};
(47.5,10)*{=F_G};
\end{xy}
\end{equation}

\vspace{.5in}

\noindent
Moreover, since $P$ is isomorphic to ${\underline{h}}^*(EG)$,
there are bundle isomorphisms
\[
E = P \times_G F \simeq {\underline{h}}^*(EG) \times_V F \simeq \underline{h}^*(EG \times_G F) = \underline{h}^*F_G.
\]
Thus, $F_G$ can be viewed as a universal fiber bundle with fiber $F$ and
structure $G$ from which all fibers bundles with fiber $F$ and structure group $G$
can be pulled back.

\section{Two Main Ideas}

In this section we isolated the two main ideas for evaluating the Gysin map $f_*\colon H^*(E) \to H^*(M)$.
Recall that a $G$-space $F$ is said to be \term{equivariantly formal} if
the canonical map $H_G^*(F) \to H^*(F)$ is surjective.
If $F$ is equivariantly formal, then by the Leray--Hirsch theorem the cohomology classes on $E$ are
 generated by pullbacks $f^*a$ of classes
$a$ from $M$ (``basic classes'') and pullbacks $h^*b$ of classes $b$ from the universal
fiber bundle $F_G$ (``fiber classes").
By the push-pull formula (\cite[Prop.~8.3]{borel--hirzebruch} or \cite[Lem.~1.5]{chern}), 
we get from \eqref{e:universal} the commutative diagram
\[
\xymatrix{
H^*(E)  \ar[d]_-{f_*}& H_G^*(F) \ar[d]^-{\pi_{G*}}\ar[l]_-{h^*} \\
H^*(M) & H^*(BG). \ar[l]^-{\underline{h}^*} 
}
\]
In other words,
\[
f_*\big( (f^*a)h^*b\big) = a f_*(h^*b) = a \underline{h}^*(\pi_{G*} b).
\]
Thus, it is enough to know how
to compute $\pi_{G*}$.  This is the first main idea.

The second main idea is based on the fact that for any $G$-space $X$,
where $G$ is a compact Lie group with maximal torus $T$ and Weyl
group $W$, we have
\[
H_G^*(X) = H_T^*(X)^W.
\]
Thus, the $G$-equivariant cohomology of $X$ injects into the $T$-equivariant comology
of $X$. 
Although the equivariant
localization theorem is valid only for a torus action,
we can apply it to get $\pi_{G*}\colon H_G^*(F) \to H^*(BG)$
for any compact, connected Lie group $G$.
To see this, first note that since 
\[
F_T = (ET \times F)/T = (EG \times F)/T
\]
and $F_G = (EG \times F)/G$, there is a natural
projection map $F_T \to F_G$ with $G/T$ as fiber.
This map fits into a commutative diagram
\[
\xymatrix{
F_G \ar[d]_-{\pi_G} & F_T \ar[l] \ar[d]^-{\pi_T}\\
\pt_G & {\ } \pt_T. \ar[l] }
\]
The push-pull formula  gives
\[
\xymatrix{
H_G^*(F) \ar[d]_-{\pi_{G^*}} \ar@{^{(}->}[r] & H_T^*(F) \ar[d]^-{\pi_{T*}} \\
H^*(BG) \ar@{^{(}->}[r] & H^*(BT). }
\]
By \cite[Lemma 4]{tu10} both horizontal maps are inclusions.
The equivariant localization formula describes the map $\pi_{T*}$ as
a finite sum.
By the commutativity of the diagram, the same is true of $\pi_{G*}$.

\section{Gysin Formula for a Complete Flag Bundle}

For a fiber bundle with equivariantly formal fiber, the method outlined
above produces a formula for the Gysin homomorphism.
As an example, consider the complete flag bundle $f\colon \Fl(V) \to
M$
associated to a complex vector bundle $V \to M$ of rank $\ell$
over a manifold $M$.

Put a Hermitian metric on $V$ and let $P$ be the principal bundle of 
unitary frames of $V$.
Let $G$ be the unitary group $\Unit(\ell)$ and $T = \Unit(1) \times \cdots
\times U(1)$ ($\ell$ times), a maximal torus in $G$.
Then $G/T = \Unit(\ell)/\big( U(1) \times \cdots \times U(1)\big)$ 
is a complete flag manifold and $f\colon \Fl(V) = P \times_G (G/T) \to M$ is the
associated complete flag bundle with fiber $G/T$.

By \eqref{e:universal}, there is a commutative diagram

\[ 
\begin{xy}
(6,5)*{\xymatrix{
\Fl(V) \ar[d]_-f \ar[r]^-h & 
EG \times_G (G/T) \ar[d]^-{\pi_G} \\
M \ar[r]_-{\underline{h}} & BG.}
};
(58,10)*{=(G/T)_G};
\end{xy}
\]

\vspace{.45in}
\noindent
The upper right corner is
\[
EG \times_G (G/T) \simeq EG/T = ET/T = BT.
\]
The push-pull formula then gives the commutative diagram
of cohomology groups

\[ 
\begin{xy}
(5,10)*{H^*(G/T)};
(6,5)*{\xymatrix{
&H^*\big(\Fl(V)\big) \ar[d]_-{f_*} \ar[l] & 
H_G^*(G/T) \ar[l]_-{h^*} \ar[d]^-{\pi_{G*}} \\
&H^*(M) & H^*(BG). \ar[l]^-{{\underline{h}}^*}
}};
(76,10)*{=H^*(BT)};
\end{xy}
\]

\vspace{.5in}
\noindent
Because $G/T$ is equivariantly formal, 
the composite map in the top line of the diagram above is surjective
and so there are global classes
on $H^*\big( \Fl(V)\big)$ that restrict to a basis on $H^*(G/T)$.
By the Leray--Hirsch theorem, the cohomology of $\Fl(V)$
is generated as an $H^*(M)$-module by the fiber classes $h^*\big( b(u)
\big)$
for $b(u) \in H^*(BT) = \R[u_1, \ldots, u_{\ell}]$.
Since $h^*$ is a ring homomorphism,
\[
h^*\big( b(u)\big) \coloneqq h^*\big( b(u_1, \ldots, u_{\ell}) \big)
= b(h^*u_1, \ldots, h^*u_{\ell})  = b(a_1, \ldots, a_{\ell}),
\]
where we write $h^*u_i=a_i$.
By the projection formula,
\[
f_*\big( (f^*c) b(a) \big) = c f_*\big(b(a)\big)  \text{ for } c \in H^*(M).
\]
Hence, the Gysin map for $f\colon \Fl(V) \to M$ is completely
determined by $f_*\big( b(a)\big)$ for fiber classes $b(a) \in H^*\big(
\Fl(V)\big)$.
Since $f^*\colon H^*(M) \to H^*\big( \Fl(V)\big)$ is injective,
$f^* f_* \big(b(a)\big)$ determines $f_*\big( b(a)\big)$.
In \cite[Prop.~13]{tu13}, we obtain the following formula for the
Gysin map of the associated complete flag bundle.

\begin{theorem} 
For the associated complete flag bundle $f\colon \Fl(V) \to M$
of a vector bundle $V \to M$,
if $b(u) \in H^*(BT)= \R[u_1, \ldots, u_{\ell}]$
and $a_i = h^*u_i$, then $b(a) \in H^*\big( \Fl(V) \big)$ and
\[
f^* f_* b(a) = \sum_{w\in S_n} w\cdot
\left(
\dfrac{b(a)}
{\prod_{i < j} (a_i - a_j)} \right),
\]
where $S_n$ is the symmetric group on $n$ letters and
$w\cdot b(a_1, \ldots, a_n) = b\big(a_{w(1)}, \ldots, a_{w(n)}\big)$.
\end{theorem}

\end{document}